\documentclass[reqno,11pt]{amsart}

\usepackage{amsmath}
\usepackage{amsfonts}   
\usepackage{amssymb}    
\usepackage{verbatim}
\usepackage{color, graphicx, graphics, pslatex, subfigure}
\usepackage{setspace}

\theoremstyle{plain}
\newtheorem{lemma}{Lemma}
\newtheorem{theorem}[]{Theorem}

\newtheorem{proposition}[]{Proposition}

\theoremstyle{remark}

\newcommand*  {\R} {{\mathbb R}}
\newcommand*  {\Z} {{\mathbb Z}}

\newcommand*{\norm}[3][{\vphantom 1}]{\lVert #2 \rVert_{#3}^{#1}}

\newcommand*{\abs}[1]{\lvert #1 \rvert}

\newcommand*{\bbabs}[1]{\big \lvert #1 \big \rvert}

\newcommand*{\ang}[1]{\left\langle #1 \right\rangle}

\newcommand{\nm}   [1]{\norm{#1}{}{}}
\newcommand{\ddt}     {\frac{d}{dt}}

\newcommand{\half}    {\frac{1}{2}}

\begin{document}

\title[Gevrey regularity of the global attractor of the $3$D Navier-Stokes-Voight equations]%
{Gevrey regularity of the global attractor of the $3$D
Navier-Stokes-Voight equations}

\date{September 20, 2007}



\author[V. K. Kalantarov]{Varga K. Kalantarov}
\address[V. K. Kalantarov]%
{Department of Mathematics, Koc University \\
 Rumelifeneri Yolu, Sariyer 34450 \\
 Istanbul, Turkey }
\email{vkalantarov@ku.edu.tr}

\author[B. Levant]{Boris Levant}
\address[B. Levant]%
{Department of Computer Science and Applied Mathematics \\
 Weizmann Institute of Science \\
 Rehovot, 76100 \\
 Israel}
\email{boris.levant@weizmann.ac.il}

\author[E. S. Titi]{Edriss S. Titi}
\address[E. S. Titi]%
{Department of Mathematics and
 Department of Mechanical and Aerospace Engineering \\
 University of California \\
 Irvine, CA 92697 \\
 USA \\
 Also, Department of Computer Science and Applied Mathematics \\
 Weizmann Institute of Science \\
 Rehovot, 76100 \\
 Israel}
\email{etiti@math.uci.edu and edriss.titi@weizmann.ac.il}

\begin{abstract}
Recently, the Navier-Stokes-Voight (NSV) model of viscoelastic
incompressible fluid has been proposed as a regularization of the
$3$D Navier-Stokes equations for the purpose of direct numerical
simulations. In this work we prove that the global attractor of the
$3$D NSV equations, driven by an analytic forcing, consists of
analytic functions. A consequence of this result is that the
spectrum of the solutions of the $3$D NSV system, lying on the
global attractor, have exponentially decaying tail, despite the fact
that the equations behave like a damped hyperbolic system, rather
than the parabolic one. This result provides an additional evidence
that the $3$D NSV with the small regularization parameter enjoys
similar statistical properties as the $3$D Navier-Stokes equations.
Finally, we calculate a lower bound for the exponential decaying
scale -- the scale at which the spectrum of the solution start to
decay exponentially, and establish a similar bound for the steady
state solutions of the $3$D NSV and $3$D Navier-Stokes equations.
Our estimate coincides with similar available lower bound for the
smallest dissipation length scale of solutions of the $3$D
Navier-Stokes equations.
\end{abstract}

\maketitle

\medskip

\medskip

\section{Introduction}

We consider the Navier-Stokes-Voight (NSV) model of viscoelastic
fluid which is governed by the system of equations
\begin{subequations}
\label{eq_nsv}
\begin{gather}
u_t - \nu \triangle u - \alpha^2 \triangle u_t + (u \cdot \nabla)
u + \nabla p = f, \label{eq_nsv_eq} \\
\nabla \cdot u = 0, \label{eq_nsv_div} \\
u(x, 0) = u^{in}(x), \label{eq_nsv_initial}
\end{gather}
\end{subequations}
in $\Omega = [0, L]^3\subset \R^3$, equipped with the periodic
boundary conditions. $u(x, t)$ represents the velocity field, $p$ is
the pressure, $\nu > 0$ stands for kinematic viscosity, $f$ is the
forcing, and finally, $\alpha$ is a real positive length scale
parameter, for which the ratio $\frac{\alpha^2}{\nu}$ characterizes
the response time that is required for the fluid to respond to the
applied force. The system (\ref{eq_nsv}) was first studied by
Oskolkov, who introduced the NSV equations (see
\cite{Oskolkov73_NSV}, \cite{Oskolkov80_Voight}) as a model of
motion of linear, viscoelastic fluid.

Recently, in \cite{CLT}, the 3D Navier-Stokes-Voight equations were
suggested as a regularization model for the 3D Navier-Stokes
equations, where $\alpha$ is considered a small regularization
parameter. First, it was recognized that the inviscid ($\nu = 0$)
version of the NSV system (\ref{eq_nsv}) coincides with the inviscid
simplified Bardina sub-grid scale model of turbulence. The viscous
simplified Bardina model was introduced and studied in \cite{L1}
(see also \cite{BFR}, and \cite{L2} for the original Bardina model).
In \cite{CLT} the viscous and inviscid simplified Bardina models
were shown to be globally well-posed. Moreover, it was also shown
that the viscous simplified Bardina model has a finite dimensional
global attractor, and the energy spectrum was investigated in
\cite{CLT}. Viewed from the numerical analysis point of view the
authors of \cite{CLT} proposed the inviscid simplified Bardina model
(or equivalently the inviscid NSV equations) as an inviscid
regularization (because no additional viscosity or hyperviscosity is
introduced) of the 3D Euler equations, subject to periodic boundary
conditions. Motivated by this observation the system (\ref{eq_nsv})
was also proposed in \cite{CLT} as a regularization, of the 3D
Navier-Stokes equations for the purpose of direct numerical
simulations for both the periodic and the no-slip Dirichlet boundary
conditions.

The addition of the $-\alpha^2 \triangle u_t$ term has two main
effects. First, it regularizes the equation in a way that the
three-dimensional system (\ref{eq_nsv}) becomes now globally
well-posed (see \cite{CLT}, \cite{Oskolkov73_NSV}). On the other
hand, as was noted in \cite{KaTi07_NSV}, it changes the parabolic
character of the original Navier-Stokes equations. Therefore, one
does not observe any immediate smoothing of the solutions, as
expected in parabolic PDEs. We also remark that this type of
inviscid regularization has been recently used for the
two-dimensional surface quasi-geostropic model
\cite{KhouiderTiti_Quasi}. In particular, necessary and sufficient
conditions for the formation of singularity were presented in terms
of regularizing parameter.

The long-time dynamics of the system (\ref{eq_nsv}) has been studied
in \cite{Ka1} and \cite{KaTi07_NSV}, where the existence of the
finite-dimensional global attractor of the system has been
established. Moreover, upper bounds for the number of determining
modes, and the fractal dimension of the global attractor of the 3D
NSV model where derived in \cite{KaTi07_NSV}. In particular, it was
shown that the attractor lies in the bounded subset of the Sobolev
space $H^1(\Omega)$, whenever the forcing term $f\in L^2(\Omega)$.

In this work we show that the global attractor of the 3D NSV model
consists of the real analytic functions, whenever the forcing term
$f$ is analytic. The idea is to construct an asymptotic
approximation $v(x, t)$ to the solution $u(x, t)$ of the system
(\ref{eq_nsv}) satisfying
\[
\lim_{t\to \infty} \nm{v(\cdot, t) - u(\cdot, t)}_{L^2(\Omega)} = 0,
\]
and show that $v(x, t)$ lies in certain Gevrey class -- a subspace
of the real analytic functions. Functions belonging to Gevrey
regularity class are characterized by the exponential decay of the
tail of their Fourier coefficients. Our method of the proof --
splitting of $v(x, t)$ into higher and lower Fourier components, has
been used before in the context of the weakly damped driven
nonlinear Schr\"{o}dinger equation in \cite{OliTi98_Schrodinger} and
a model of B\'{e}nard convection in a porous medium in
\cite{OliTi00_Porous} (see also \cite{Goubet96_Schrodinger}).
Recently, the authors of \cite{CPS04_Gevrey} followed the same
methods to prove the Gevrey regularity of the global attractor of
the generalized Benjamin-Bona-Mahony equation.

An important consequence of our result is that the solutions of the
3D NSV system (\ref{eq_nsv}) lying on the global attractor posses a
dissipation range, despite the fact that the equations behave like
the damped hyperbolic system, rather than the parabolic equation.
This fact provides an additional evidence that (\ref{eq_nsv}), with
the small regularization parameter $\alpha$, can indeed be used as a
model to study the statistical properties of turbulent solutions of
the 3D Navier-Stokes equations, a subject of ongoing research.

Finally, we obtain bounds for the exponential decaying length scale,
which is related to the dissipation length scale,  of the general
solutions of the NSV system lying on the global attractor. The
obtained estimate is similar to the bounds for the smallest length
scale in the turbulent flow that was previously calculated for the
solutions of the $3$D Navier-Stokes equations in \cite{DT95_Decay}.
In addition, using the techniques introduced in \cite{OliTi01_Stat},
we estimate the exponential decaying scale of the stationary
solutions of the $3$D NSV and $3$D Navier-Stokes equations. Our
bounds coincide with those obtained in this paper for the general
solutions of the NSV system lying on the global attractor, and for
those of the $3$D Navier-Stokes equations reported in
\cite{DT95_Decay}.

\section{Preliminaries}

In this paper we will use the following notations, which are
standard in the mathematical theory of the Navier-Stokes equations
(see, e.g., \cite{CF88}, \cite{FMRT01}, \cite{Te01_NSE}).

Let $\Omega = [0, L]^3$. We denote by $L^p(\Omega)$, for $1 \le p
\le \infty$, and $H^m(\Omega)$ -- the usual Lebesgue and Sobolev
spaces of the periodic functions on $\Omega$ respectively. Let
$\mathcal{F}$ be the set of all vector trigonometric polynomials on
the periodic domain $\Omega$, and denote
\begin{equation*}
\mathcal{V} = \Big \{ \varphi\in \mathcal{F} \;:\; \nabla\cdot
\varphi = 0, \; \text{and} \; \int_\Omega \varphi(x) dx = 0 \Big \}.
\end{equation*}
We set $H$, and $V$ to be the closures of $\mathcal{V}$ in the
$L^2(\Omega)$ and $H^1(\Omega)$ topology respectively.

We denote by $P_\sigma : L^2 \to H$ -- the Helmholtz-Leray
orthogonal projection operator, and by $A = - P_\sigma \triangle$ --
the Stokes operator subject to the periodic boundary conditions with
domain $D(A) = (H^2(\Omega))^3\cap V$. Observe, that in the
space-periodic case
\[
A u = - P_\sigma \triangle u = - \triangle u, \;\;\; \text{for all}
\; u\in D(A).
\]
The operator $A^{-1}$ is a positive definite, self-adjoint, compact
operator from $H$ into $H$. We denote by $0 < \Big( \frac{2\pi}{L}
\Big)^3 = \lambda_1 \le \lambda_2 \le \dots$ the eigenvalues of $A$,
repeated according to their multiplicities. The eigenvalues
$\lambda_j$ satisfy, for some dimensionless constant $c_0 > 0$,
\[
\frac{j^{2/3}}{c_0} \le \frac{\lambda_j}{\lambda_1} \le c_0 j^{2/3},
\;\;\; \text{for} \; j = 1, 2, 3, \dots.
\]
In the periodic case this observation is simple (see, e.g.,
\cite{CF88}), however, in the general case, this is a result of the
famous Weyl's formula for the case of the Stokes operator due to
M\'{e}tivier (see, e.g., \cite{CF88}, \cite{Metiv78}, \cite{Te88}).

For any $s\in \R$, we can define the Hilbert spaces $V_s :=
D(A^{s/2})$ with the inner product and norm
\begin{equation*}
(u, v)_s = \sum_{j\in \Z^3}  u_j\cdot v_j \abs{j}^{2s}, \;\;\;
\abs{u}_s^2 = (u, u)_s,
\end{equation*}
for every $u, v\in V_s$, where $u_j, v_j$ are the corresponding
Fourier coefficients of $u$ and $v$ respectively. Note, that $V_0 =
H$. We will denote the corresponding inner product and norm in $H$
by $(\cdot,\cdot)$ and $|\cdot|$, respectively.  Moreover, we denote
$V = V_1$, and the corresponding inner product and norm will be
written for $u, v \in V$
\[
((u, v)) = (u, v)_1, \;\;\; \nm{u} = \abs{u}_1.
\]

For any $w_1, w_2\in \mathcal{V}$ we define the following bilinear
form
\begin{equation*}
B(w_1, w_2) = P_\sigma \big((w_1\cdot \nabla) w_2 \big).
\end{equation*}
It can be shown (see, e.g., \cite{CF88}, \cite{Te01_NSE}) that $B$
can be extended to a continuous map $B: V\times V\to V'$, where $V'
= V_{-1}$ is a dual space of $V$. In particular, for $u, v, w\in V$,
there exists a constant $c > 0$, depending only on $\Omega$, such
that
\begin{equation} \label{eq_nlin}
\bbabs{ \ang{ B(u, v), w}_{V'} } \le c \lambda_1^{-3/4}
\abs{u}^{1/2} \nm{u}^{1/2} \nm{v} \nm{w},
\end{equation}
where $\ang{x, y}_{V'}$ denotes an action of an element $x\in V$ on
the element of the dual space $y\in V'$.

Finally, using the above definitions, we write the system
(\ref{eq_nsv}) in the following equivalent functional form
\begin{subequations}
\label{eq_nsvabs}
\begin{gather}
u_t + \nu A u + \alpha^2 A  u_t + B(u, u) = f, \label{eq_nsvabs_eq} \\
u(x, 0) = u^{in}(x). \label{eq_nsvabs_initial}
\end{gather}
\end{subequations}

To show that the solution of the problem (\ref{eq_nsvabs}) has an
analytic asymptotic (in time) approximation, we will use the concept
of the Gevrey class regularity. For a given $\tau
> 0$, and $r\ge 0$, we define the Gevrey space to be
\begin{equation*}
G_\tau^{r} := D(A^{r/2} e^{\tau A^{1/2}}) = \{ u\in H \;:\; \bbabs{
A^{r/2} e^{\tau A^{1/2}} u}^2 = \sum_{j\in \Z^3} \abs{u_j}^2
\abs{j}^{2r} e^{2 \tau \abs{j}} < \infty \}.
\end{equation*}
The space is equipped with the corresponding inner product and norm
\begin{equation*}
(u, v)_{r, \tau} = (A^{r/2} e^{\tau A^{1/2}} u, A^{r/2} e^{\tau
A^{1/2}} v) = \sum_{j\in \Z^3} u_j\cdot v_j \abs{j}^{2r} e^{2 \tau
\abs{j}}, \;\;\; \abs{u}_{r, \tau} = \abs{A^{r/2} e^{\tau A^{1/2}}
u},
\end{equation*}
for $u, v\in G_\tau^{r}$. One can prove that the space of real
analytic functions $C^\omega(\Omega)$ has the following
characterization
\begin{equation*}
C^\omega(\Omega) = \bigcup_{\tau > 0} G_\tau^{r},
\end{equation*}
for any $r \ge 0$ (see, e.g., \cite{LevOliver97_Gevrey}). The
concept of the Gevrey class regularity for showing the analyticity
of the solutions of the Navier-Stokes equations, was first
introduced in \cite{FT89_Gevrey}, simplifying earlier proofs. Later
this technique was extended to the large class of analytic nonlinear
parabolic equations in \cite{FT98_Gevrey}.

We conclude this section by a few technical propositions that will
be used in the proof of our main results. First, we will need the
following estimates for the nonlinear term. The proof of Proposition
\ref{prop_nonlin}, below, is achieved by standard interpolation
estimates using the Gagliardo-Nirenberg and Ladyzhenskaya
inequalities (see, e.g., \cite{CF88}, \cite{Te01_NSE}).

\begin{proposition} \label{prop_nonlin}
The bilinear form $B(u, u)$ satisfies:
\begin{enumerate}
\item If $u\in V$, then $B(u, u)\in V_{-1/2}$, and
\begin{equation} \label{eq_nl1}
\bbabs{B(u, u)}_{-1/2} \le c_1 \lambda_1^{-3/4} \nm{u}^2.
\end{equation}

\item If $u\in V_{3/2}$, then $B(u, u)\in H$, and
\begin{equation} \label{eq_nl2}
\bbabs{B(u, u)} \le c_2 \lambda_1^{-3/4} \nm{u} \abs{u}_{3/2}.
\end{equation}

\item For any integer $m\ge 1$, if $u\in V_{m+1}$, then $B(u, u)\in V_{m}$,
and
\begin{equation} \label{eq_nl3}
\bbabs{B(u, u)}_{m} \le c_m \lambda_1^{-7/8} \nm{u}^{1/4}
\abs{u}_{2}^{3/4} \abs{u}_{m+1},
\end{equation}
\end{enumerate}
where $c_1, c_2, c_m > 0$ are scale invariant constants, and $c_m$
depends on $m$.
\end{proposition}

Let $\lambda > 0$, denote by $P_\lambda$ the $H$-orthogonal
projection onto the span of eigenfunctions of $A$ corresponding to
eigenvalues of the magnitude less then or equal to $\lambda$. Denote
$Q_\lambda = I - P_\lambda$. The following Poincar\'{e}-type
inequalities hold.

\begin{proposition} \label{prop_poincare}
Let $\bar v\in P_\lambda G_\tau^{r+1}$, and $\hat v\in Q_\lambda
G_\tau^{r+1}$. Then,
\begin{equation} \label{eq_poincare}
\abs{\bar v}_{r+1, \tau} \le e^{\tau \lambda^{1/2}} \abs{\bar
v}_{r+1}, \;\;\; \text{and} \;\;\; \abs{\hat v}_{r, \tau} \le
\lambda^{-1/2} \abs{\hat v}_{r+1, \tau}.
\end{equation}
\end{proposition}

We will also need an estimate for the nonlinear term in the Gevrey
space. Similar inequalities can be found in \cite{FMRT01} (see also
\cite{DT95_Decay}, \cite{FT89_Gevrey}).

\begin{proposition} \label{prop_gevreyineq}
For any $\tau > 0$, $u, w\in G_\tau^{2}$, and $v\in G_\tau^{1}$, the
following inequality holds
\begin{equation} \label{eq_gevreyineq}
\bbabs{\big (B(u, v), w \big)_{1, \tau} } \le C_1 \lambda_1^{-3/4}
\abs{u}_{1, \tau}^{1/2} \abs{u}_{2, \tau}^{1/2} \abs{v}_{1, \tau}
\abs{w}_{2, \tau},
\end{equation}
for some scale invariant constant $C_1 > 0$.

\end{proposition}

It is not difficult to prove the following proposition using the
Galerkin approximation procedure (see \cite{KaTi07_NSV}).

\begin{proposition} \label{prop_solution}
Let $s\in \R$. Assume that $g(t)\in L^\infty([0, T], V_{s-2})$, for
some $0 < T < \infty$. Then the linear problem
\begin{equation*}
z_t + \nu A z + \alpha^2 A z_t = g(t), \;\;\; z(0) = 0,
\end{equation*}
has a unique solution $z(t)\in C([0, T], V_s)$. In addition, the
following estimate holds,
\begin{equation} \label{eq_zbnd}
\abs{z(t)}_s \le \frac{\nm{g(t)}_{L^\infty([0, T], V_{s-2})}}{\alpha
\nu \sqrt{d_0}},
\end{equation}
for all $t\in [0, T]$, and $d_0 = \big( \frac{1}{\lambda_1} +
\alpha^2 \big)^{-1}$.
\end{proposition}

We will use the following proposition which we state here without a
proof.

\begin{proposition} \label{prop_bnd}
Let  $\varphi(t)$ be a nonnegative absolutely continuous function on
$[t_0, \infty)$, for some $t_0 \ge 0$, satisfying, for all $t\ge
t_0$, the  inequality
\[
\frac{d \varphi}{dt} \le -a \varphi + b \varphi^{3/2} + c, \;\;\;
\varphi(t_0) = 0.
\]
Where the positive constant coefficients $a, b, c$ obey the
inequality
\begin{equation} \label{eq_coeff}
b c^{1/2} < \Big( \frac{a}{2} \Big)^{3/2}.
\end{equation}
Then for all $t\ge t_0$
\[
\varphi(t)\le \frac{2 c}{a}.
\]
\end{proposition}

Finally, we will need the following Lemma from \cite{JT92_2dnse}
(see also \cite{FMRT01}).

\begin{lemma} \label{lem_1}
Let $a(t)$ and $b(t)$ be locally integrable functions on
$(0,\infty)$ which satisfy for some $T > 0$ the conditions
\begin{align*}
\liminf_{t\rightarrow \infty} \frac{1}{T} \int_{t}^{t+T} a(\tau)
d\tau & > 0,
\\
\limsup_{t\rightarrow \infty} \frac{1}{T} \int_{t}^{t+T} a^{-}(\tau)
d\tau & < \infty,
\\
\limsup_{t\rightarrow \infty} \frac{1}{T} \int_{t}^{t+T} b^{+}(\tau)
d\tau & = 0,
\end{align*}
where $a^{-}=\max\{-a,0\}$ and $b^{+}=\max\{b,0\}$. Suppose that
$\phi(t)$ is a nonnegative, absolutely continuous function on $[0,
\infty)$ that satisfies the following inequality, almost everywhere
on $[0, \infty)$,
\[
\phi'(t) + a(t) \phi(t) \leq b(t).
\]
Then $\phi(t)\to 0$, as $t\to \infty$.
\end{lemma}

\section{Asymptotic approximation in $V_m$}

The question of global existence and uniqueness of \eqref{eq_nsv}
was first studied  in \cite{Oskolkov73_NSV} (the inviscid case, $\nu
= 0$, was studied in \cite{CLT}). It was shown that for every
$u^{in}\in V$, the solution of the the system (\ref{eq_nsv}) is
globally well-posed, and satisfies $u(x, t)\in L^\infty([0, \infty),
V)$. In this section we construct an asymptotic approximation of the
solution of (\ref{eq_nsv}) in the space $V_m$, for every $m\ge 2$.
The result can be stated as follows.

\begin{theorem} \label{thm_vm}
Let $m \ge 2$ be an integer. Consider a solution $u(x, t)$ of the
NSV system (\ref{eq_nsv}), corresponding to the initial condition
$u^{in}\in V$ with the forcing $f\in V_{m-2}$. Then there exists a
function
\begin{equation*}
v^{(m)}(t)\in L^\infty([0, \infty), V_m),
\end{equation*}
satisfying
\begin{equation*}
\lim_{t\to \infty} \nm{u(t) - v^{(m)}(t)} = 0.
\end{equation*}
\end{theorem}

\begin{proof}
Let us fix $m\ge 2$, and let $u^{in}\in V$. First, let us write the
solution $u(t) = v(t) + w(t)$, where $v(t)$ and $w(t)$ satisfy the
coupled system
\begin{subequations}
\label{eq_approx1}
\begin{gather}
v_t + \nu A v + \alpha^2 A v_t = f - B(u, u), \;\;\;\;\; v(0) = 0,
\label{eq_approx1v} \\ w_t + \nu A w + \alpha^2 A w_t = 0,
\;\;\;\;\; w(0) = u^{in}. \label{eq_approx1w}
\end{gather}
\end{subequations}
This decomposition has been used in \cite{KaTi07_NSV}. First, by
using the fact that $u(x, t)\in L^\infty([0, \infty), V)$, and
applying subsequently the first part of
Proposition~\ref{prop_nonlin} and Proposition~\ref{prop_solution} to
equation (\ref{eq_approx1v}), we conclude that
\begin{equation} \label{eq_v}
v(t)\in L^\infty([0, \infty), V_{3/2}).
\end{equation}
Next, from equation (\ref{eq_approx1w}) we immediately get
\begin{equation} \label{eq_w}
\abs{w(t)}^2 + \alpha^2 \nm{w(t)}^2 \le e^{-\nu d_0 t} \big(
\abs{u^{in}}^2 + \alpha^2 \nm{u^{in}}^2 \big),
\end{equation}
where $d_0 = \big( \frac{1}{\lambda_1} + \alpha^2 \big)^{-1}$.
Therefore, $v(x, t)$ is an asymptotic (in time) approximation of
$u(x, t)$, namely
\[
\lim_{t\to \infty} \nm{u(t) - v(t)} = \lim_{t\to \infty} \nm{w(t)} =
0.
\]

At the next step, let us consider $v^{(2)}(x, t)$ -- the solution of
the following equation
\begin{equation} \label{eq_approx2}
v^{(2)}_t + \nu A v^{(2)} + \alpha^2 A v^{(2)}_t = f - B(v, v),
\;\;\;\;\; v^{(2)}(0) = 0.
\end{equation}
According to the second part of Proposition~\ref{prop_nonlin}, the
right-hand side of equation (\ref{eq_approx2}) is in $L^\infty([0,
\infty), H)$. Therefore, applying Proposition~\ref{prop_solution} we
conclude that the unique solution of equation (\ref{eq_approx2})
satisfies
\begin{equation} \label{eq_v'}
v^{(2)}(t)\in L^\infty([0, \infty), V_2).
\end{equation}
Denote $z^{(2)} = v^{(2)} - v$, which satisfies
\begin{equation} \label{eq_asym2}
z^{(2)}_t + \nu A z^{(2)} + \alpha^2 A z^{(2)}_t = B(u, u - v) + B(u
- v, v), \;\;\;\;\; z^{(2)}(0) = 0.
\end{equation}
According to Proposition~\ref{prop_solution} equation
(\ref{eq_asym2}) has a unique solution $z^{(2)}(t)\in L^\infty([0,
\infty), V_{3/2})$. This is because $u\in L^\infty([0, \infty), V)$,
and $v$ satisfies (\ref{eq_v}). Therefore, we can take an inner
product of equation (\ref{eq_asym2}) with $z$. Using inequality
(\ref{eq_nlin}) we get
\begin{align}
\half \ddt \big( \abs{z^{(2)}(t)}^2 & + \alpha^2 \nm{z^{(2)}(t)}^2
\big) + \nu \nm{z^{(2)}(t)}^2 = \notag
\\
& = \big(B(u(t), u(t) - v(t)), z^{(2)}(t) \big) + \big(B(u(t) -
v(t), v(t)), z^{(2)}(t)\big) \le \notag
\\
& \le c \nm{u(t) - v(t)} \big( \nm{u(t)} + \nm{v(t)} \big)
\nm{z^{(2)}(t)} \le \notag
\\
& \le \frac{c^2}{2 \nu} \nm{u(t) - v(t)}^2 \big( \nm{u(t)} +
\nm{v(t)} \big)^2 + \frac{\nu}{2} \nm{z^{(2)}(t)}^2,
\label{eq_asym21}
\end{align}
where the last relation follows from Young's inequality. Finally, we
get
\begin{align*}
\ddt \big( \abs{z^{(2)}(t)}^2 + \alpha^2 \nm{z^{(2)}(t)}^2 \big) & +
\frac{\nu d_0}{2} \big( \abs{z^{(2)}(t)}^2 + \alpha^2
\nm{z^{(2)}(t)}^2 \big) \le
\\
& \le \frac{c^2}{2 \nu} \nm{u(t) - v(t)}^2 \big( \nm{u(t)} +
\nm{v(t)} \big)^2.
\end{align*}
Using the fact that $u(t), v(t)$ are bounded uniformly in time in
the $V$ norm, and that
\[
\nm{u(t) - v(t)} = \nm{w(t)} \to 0, \;\;\; \text{as} \;\;\; t\to
\infty,
\]
we conclude, after applying Lemma~\ref{lem_1} that
\begin{equation*}
\lim_{t\to \infty} \nm{z^{(2)}(t)} = \lim_{t\to 0} \nm{v(t) -
v^{(2)}(t)} = \lim_{t\to 0} \nm{u(t) - v^{(2)}(t)} = 0. \notag
\end{equation*}

We can continue by induction. Fix $2 \le n \le m$, and assume that
we have constructed $v^{(j)}(t)\in L^\infty([0, \infty), V_{j})$,
for $j = 2, 3, \dots, n-1$, such that for any $j$
\begin{equation} \label{eq_j}
\lim_{t\to \infty} \nm{v^{(j-1)}(t) - v^{(j)}(t)} = \lim_{t\to 0}
\nm{u(t) - v^{(j)}(t)} = 0.
\end{equation}
Let consider the following equation
\begin{equation} \label{eq_approxm}
v^{(n)}_t + \nu A v^{(n)} + \alpha^2 A v^{(n)}_t = f - B(v^{(n-1)},
v^{(n-1)}), \;\;\;\;\; v^{(n)}(0) = 0.
\end{equation}
Then, according to Proposition~\ref{prop_solution}, and due to the
estimates on the nonlinear term of Proposition~\ref{prop_nonlin},
the unique solution $v^{(n)}(t)$ of the equation (\ref{eq_approxm})
satisfies $v^{(n)}(t)\in L^\infty([0, \infty), V_n)$. Moreover,
denote $z^{(n)} = v^{(n)} - v^{(n-1)}$, satisfying
\begin{align*}
z^{(n)}_t + \nu A z^{(n)} + \alpha^2 A z^{(n)}_t & = B(v^{(n-2)},
v^{(n-2)} - v^{(n-1)}) + B(v^{(n-2)} - v^{(n-1)}, v^{(n-1)}),
\\
z^{(n)}(0) & = 0.
\end{align*}
Taking the inner product of the last equation with $z^{(n)}(t)$ and
using Proposition~\ref{prop_nonlin} and relation (\ref{eq_j}) we can
show by Lemma~\ref{lem_1} that
\begin{equation*}
\lim_{t\to \infty} \nm{z^{(n)}(t)} = \lim_{t\to 0} \nm{v^{(n-2)}(t)
- v^{(n-1)}(t)} = \lim_{t\to 0} \nm{u(t) - v^{(n)}(t)} = 0,
\end{equation*}
finishing the proof of the Theorem.

\end{proof}

It can be proved (see \cite{KaTi07_NSV}) that the solution of the
NSV equations (\ref{eq_nsv}) satisfies for all $t\ge 0$
\begin{equation}
\nm{u(t)}^2 \le \frac{e^{-\nu d_0 t}}{\alpha^2} \Bigg(
\abs{u^{in}}^2 + \alpha^2 \nm{u^{in}} - \frac{\abs{f}_{-1}^2}{\nu^2
d_0} \Bigg) + \frac{\abs{f}_{-1}^2}{\alpha^2 \nu^2 d_0}.
\end{equation}
Therefore, there exists $t_0$, depending on $\abs{u^{in}},
\nm{u^{in}}, \abs{f}_{-1}, \nu, \alpha$, and $\lambda_1$, such that
for all $t\ge t_0$
\begin{equation} \label{eq_v1bnd}
\nm{u(t)} \le M_1 := \frac{2 \abs{f}_{-1}}{\alpha \nu \sqrt{d_0}}.
\end{equation}
The following Lemma gives similar bounds for the asymptotic (in
time) approximations $v^{(m)}(x, t)$ in the corresponding norms.

\begin{lemma} \label{lem_2}
Let $f\in V_{m-2}$. Consider $t_0 \ge 0$, such that the solution of
the NSV equations (\ref{eq_nsv}) satisfies the inequality
(\ref{eq_v1bnd}) for all $t\ge t_0$. Then the following statements
are true:
\begin{enumerate}
\item The function $v(x, t)\in L^\infty([0,
\infty), V_{3/2})$, constructed in Theorem~\ref{thm_vm}, satisfies
for all $t\ge t_0$
\begin{equation} \label{eq_v32bnd}
\abs{v(t)}_{3/2} \le M_{3/2} := \frac{1}{\alpha \nu \sqrt{d_0}}
\big( \abs{f}_{-1/2} + c_1 \lambda_1^{-3/4} M_1^2 \big).
\end{equation}

\item The function $v^{(2)}(x, t)\in L^\infty([0,
\infty), V_2)$, constructed in Theorem~\ref{thm_vm}, satisfies for
all $t\ge t_0$
\begin{equation} \label{eq_v2bnd}
\abs{v^{(2)}(t)}_2 \le M_2 := \frac{1}{\alpha \nu \sqrt{d_0}} \big(
\abs{f} + c_2 \lambda_1^{-3/4} M_1 M_{3/2} \big).
\end{equation}

\item Let $m> 2$ be an integer. The function $v^{(m)}(x, t)\in L^\infty([0,
\infty), V_{m})$, constructed in Theorem~\ref{thm_vm}, satisfies for
all $t\ge t_0$
\begin{equation} \label{eq_vmbnd}
\abs{v^{(m)}(t)}_m \le M_m := \frac{1}{\alpha \nu \sqrt{d_0}} \big(
\abs{f}_{m-2} + c_m \lambda_1^{-7/8} M_1^{1/4} M_2^{3/4} M_{m-1}
\big).
\end{equation}

\end{enumerate}
\end{lemma}

\begin{proof}
Recall that $v(t)$ satisfies equation (\ref{eq_approx1}),
$v^{(2)}(t)$ satisfies (\ref{eq_approx2}). In general, $v^{(m)}(t)$,
for $m > 2$, satisfies equation (\ref{eq_approxm}). Therefore, the
proof of the Lemma is an immediate application of
Proposition~\ref{prop_solution}, in particular relation
(\ref{eq_zbnd}), and the inequalities of
Proposition~\ref{prop_nonlin}.

\end{proof}

\section{Asymptotic approximation in the Gevrey space $G_\tau^{1}$}

The results of the previous section show that with a smooth enough
forcing the global attractor of the system (\ref{eq_nsv}) lies in
$C^\infty(\Omega)$, whenever $f$ is $C^\infty(\Omega)$. However, our
goal is to show that the global attractor is real analytic, whenever
$f$ is real analytic. For this purpose we use the idea of
\cite{OliTi98_Schrodinger} and \cite{OliTi00_Porous}, to construct
the asymptotic approximation of the solution of (\ref{eq_nsv}) in
the Gevrey class $G_\tau^{2}$, for some $\tau > 0$.

\begin{theorem} \label{thm_gv}
Let $u(x, t)$ be a solution of the NSV system (\ref{eq_nsv}),
corresponding to the initial condition $u^{in}\in V$ with the
forcing $f\in G_{\tau_0}^{1}$, for some $\tau_0 > 0$. Let $t_0\ge 0$
be as in Lemma~\ref{lem_2}, then there exists a function
\begin{equation}
v^\omega(t)\in L^\infty([t_0, \infty), G_\tau^2),
\end{equation}
for some $\tau > 0$, depending only on $\abs{f}_{1, \tau_0}$, $\nu$,
$\lambda_1$ and $\alpha$, satisfying
\begin{equation} \label{eq_asymgv}
\lim_{t\to \infty} \nm{u(t) - v^\omega(t)} = 0.
\end{equation}
\end{theorem}

\begin{proof}
Let $\lambda > 0$ to be chosen later. First, consider $v^{(2)}(x,
t)$ -- an asymptotic approximation of $u(x, t)$, which is
constructed in Theorem~\ref{thm_vm}. Moreover, according to
Lemma~\ref{lem_1}, there exists a constant $M_2 > 0$ (see relation
(\ref{eq_v2bnd})), such that
\begin{equation} \label{eq_vrbound}
\abs{v^{(2)}(t)}_{2} \le M_2, \;\;\; \forall t \ge t_0.
\end{equation}

Denote $\bar v(t) = P_\lambda v^{(2)}(t)$, and consider $\hat v(t)$
-- a solution of the following equation
\begin{equation} \label{eq_hatu}
{\hat v}_t + \nu A {\hat v} + \alpha^2 A {\hat v}_t + Q_\lambda
B(\bar v + \hat v, \bar v + \hat v) = \hat f, \;\;\;\;\; {\hat
v}(t_0) = 0,
\end{equation}
for $t\ge t_0$, where, for notation simplicity, we denoted $\hat f =
Q_\lambda f$. The equation (\ref{eq_hatu}) formally looks like a
projection of the system (\ref{eq_nsv}) onto the higher wavenumber
components, however, the low wavenumber modes $\bar v$ of the
advection term satisfy a slightly different equation (see also
\cite{OlsonTi03_Data} for such a construction for studying data
assimilation). Let us denote by
\begin{equation} \label{eq_uanalytic}
v^\omega(t) = \bar v(t) + \hat v(t),
\end{equation}
for $t\ge t_0$. Our goal is to show first that there exists $\tau >
0$ such that $v^\omega\in G_\tau^{2}$. Observe, that $\bar v$ is
just a trigonometric polynomial, and in particular, is analytic.
Therefore, we need to show that we can choose $\lambda$ large
enough, such that $\hat v\in G_\tau^{2}$, for some $\tau > 0$.
Finally, we will show that $v^\omega(x, t)$ is indeed an asymptotic
approximation of $u(x, t)$.

Note, that in order to prove that the solution of the equation
(\ref{eq_hatu}) lies in a Gevrey class of real analytic functions we
consider the Galerkin procedure to equation (\ref{eq_hatu}).
However, we omit this standard procedure, and obtain formal a-priori
estimates on the solutions in the relevant Gevrey space norm. Taking
formally the inner product of the equation (\ref{eq_hatu}) in
$G_\tau^{1}$ with $\hat v$ we obtain the following inequality
\begin{align}
\half \ddt & \Big( \abs{{\hat v}}_{1, \tau}^2 + \alpha^2 \abs{{\hat
v}}_{2, \tau}^2 \Big) + \nu \abs{{\hat v}}^2_{2, \tau} \le \notag
\\
& \le \abs{(\hat f, \hat v)_{1, \tau}} + \abs{(B(\bar v, \bar v),
\hat v)_{1, \tau}} + \abs{(B(\bar v, \hat v), \hat v)_{1, \tau}} +
\abs{(B(\hat v, \bar v), \hat v)_{1, \tau}} + \abs{(B(\hat v, \hat
v), \hat v)_{1, \tau}}. \label{eq_hatnormu}
\end{align}

Next, we estimate the terms on the right-hand side of
(\ref{eq_hatnormu}). First, using subsequently the Cauchy-Schwartz
and Young inequalities, as well as Proposition~\ref{prop_poincare}
we get, assuming $\tau \le \tau_0$,
\begin{equation} \label{eq_rel1}
\abs{(\hat f, \hat v)_{1, \tau}} \le \abs{\hat f}_{1, \tau} \cdot
\abs{\hat v}_{1, \tau} \le \frac{5}{4 \nu \lambda} \abs{\hat f}_{1,
\tau}^2 + \frac{\nu}{5} \abs{\hat v}_{2, \tau}^2.
\end{equation}

Next, using Proposition~\ref{prop_gevreyineq}, Young inequality, and
the Poincar\'{e}-type inequalities of
Proposition~\ref{prop_poincare}, we get the following series of
estimates for all $t\ge t_0$
\begin{align}
\abs{(B(\bar v, \bar v), \hat v)_{1, \tau}} & \le C_1
\lambda_1^{-3/4} \abs{\bar v}_{1, \tau}^{3/2} \abs{\bar v}_{2,
\tau}^{1/2} \abs{\hat v}_{2, \tau} \le \notag
\\
& \le \frac{5 C_1^2 \abs{\bar v}_{1, \tau}^3 \abs{\bar v}_{2,
\tau}}{4 \nu \lambda_1^{3/2}} + \frac{\nu}{5} \abs{\hat v}_{2,
\tau}^2 \le \frac{5 C_1^2 e^{4 \tau \lambda^{1/2}} M_1^3 M_2}{4 \nu
\lambda_1^{3/2}} + \frac{\nu}{5} \abs{\hat v}_{2, \tau}^2.
\label{eq_rel2}
\end{align}

\begin{equation} \label{eq_rel3}
\abs{(B(\bar v, \hat v), \hat v)_{1, \tau}} \le C_1 \lambda_1^{-3/4}
\abs{\bar v}_{1, \tau}^{1/2} \abs{\bar v}_{2, \tau}^{1/2} \abs{\hat
v}_{1, \tau} \abs{\hat v}_{2, \tau} \le \frac{C_1 e^{\tau
\lambda^{1/2}} M_1^{1/2} M_2^{1/2}}{\lambda^{1/2} \lambda_1^{3/4}}
\abs{\hat v}_{2, \tau}^2.
\end{equation}

\begin{equation} \label{eq_rel4}
\abs{(B(\hat v, \bar v), \hat v)_{1, \tau}} \le C_1 \lambda_1^{-3/4}
\abs{\hat v}_{1, \tau}^{1/2} \abs{\hat v}_{2, \tau}^{3/2} \abs{\bar
v}_{1, \tau} \le \frac{C_1 e^{\tau \lambda^{1/2}} M_1}{\lambda^{1/4}
\lambda_1^{3/4}} \abs{\hat v}_{2, \tau}^2.
\end{equation}

\begin{align}
\abs{(B(\hat v, \hat v), \hat v)_{1, \tau}} & \le
\frac{C_1}{\lambda_1^{3/4}} \abs{\hat v}_{1, \tau}^{3/2} \abs{\hat
v}^{3/2}_{2, \tau} \le \frac{C_1}{\lambda^{3/4} \lambda_1^{3/4}}
\abs{\hat v}^{3}_{2, \tau} \le \frac{C_1}{\lambda^{3/4}
\lambda_1^{3/4} \alpha^3} \Big( \abs{\hat v}^{2}_{1, \tau} +
\alpha^2 \abs{\hat v}^{2}_{2, \tau} \Big)^{3/2}. \label{eq_rel5}
\end{align}

Let us set $\tau = \min \{ \lambda^{-1/2}, \tau_0 \}$. Then, we will
choose $\lambda$ large enough satisfying
\begin{equation} \label{eq_lambda1}
\max \Bigg\{ \frac{C_1 e M_1^{1/2} M_2^{1/2}}{\lambda^{1/2}
\lambda_1^{3/4}}, \frac{C_1 e M_1}{\lambda^{1/4} \lambda_1^{3/4}}
\Bigg\} \le \frac{\nu}{5}.
\end{equation}
Using the last bounds, we are ready to substitute relations
(\ref{eq_rel1}), (\ref{eq_rel2}), (\ref{eq_rel3}), (\ref{eq_rel4}),
and (\ref{eq_rel5}) into equation (\ref{eq_hatnormu}). After
rearranging the terms, we get
\begin{align}
\half \ddt \Big( \abs{{\hat v}}_{1, \tau}^2 + \alpha^2 \abs{{\hat
v}}_{2, \tau}^2 \Big) & + \frac{\nu}{5} \abs{{\hat v}}^2_{2, \tau}
\le \notag
\\
& \le \frac{C_1}{\lambda^{3/4} \lambda_1^{3/4} \alpha^3} \Big(
\abs{\hat v}^2_{1, \tau} + \alpha^2 \abs{\hat v}^2_{2, \tau}
\Big)^{3/2} + \frac{5 \abs{\hat f}_{1, \tau}^2}{4 \nu \lambda} +
\frac{5 C_1^2 e^{4} M_1^3 M_2}{4 \nu \lambda_1^{3/2}}.
\label{eq_hatnormu2}
\end{align}
Next, using Poincar\'{e}-type inequality,
Proposition~\ref{prop_poincare}, and setting $d_2 =
(\frac{1}{\lambda} + \alpha^{2})^{-1}$, we can write
\begin{equation} \label{eq_nud2}
\frac{\nu}{5} \abs{{\hat v}}^2_{2, \tau} \ge \frac{\nu d_2}{5}
(\lambda^{-1} \abs{{\hat v}}^2_{2, \tau} + \alpha^2 \abs{{\hat
v}}^2_{2, \tau}) \ge \frac{\nu d_2}{5} (\abs{{\hat v}}^2_{1, \tau} +
\alpha^2 \abs{{\hat v}}^2_{2, \tau}).
\end{equation}
Substituting into equation (\ref{eq_hatnormu2}) gives us the
inequality
\begin{align}
\half \ddt \Big( \abs{{\hat v}}_{1, \tau}^2 + \alpha^2 \abs{{\hat
v}}_{2, \tau}^2 \Big) \le & - \frac{\nu d_2}{5} \Big(\abs{{\hat
v}}^2_{1, \tau} + \alpha^2 \abs{{\hat v}}^2_{2, \tau} \Big) + \notag
\\
& + \frac{C_1}{\lambda^{3/4} \lambda_1^{3/4} \alpha^3} \Big(
\abs{\hat v}^2_{1, \tau} + \alpha^2 \abs{\hat v}^2_{2, \tau}
\Big)^{3/2} + \frac{5 \abs{\hat f}_{1, \tau}^2}{4 \nu \lambda} +
\frac{5 C_1^2 e^{4} M_1^3 M_2}{4 \nu \lambda_1^{3/2}}.
\label{eq_hatnormu3}
\end{align}
Now we can apply Proposition~\ref{prop_bnd} to the function $\varphi
(t)=\big( \abs{{\hat v(t)}}_{1, \tau}^2 + \alpha^2 \abs{{\hat
v(t)}}_{2, \tau}^2 \big)$ which is satisfying inequality
(\ref{eq_hatnormu3}). Using relation (\ref{eq_coeff}) we conclude,
that $\big( \abs{{\hat v(t)}}_{1, \tau}^2 + \alpha^2 \abs{{\hat
v(t)}}_{2, \tau}^2 \big)$ is bounded for all $t\ge t_0$, and in
particular $v(t)\in L^\infty([t_0, \infty), G^2_\tau)$, whenever the
following holds
\begin{equation*}
\frac{C_1}{\lambda^{3/4} \lambda_1^{3/4} \alpha^3} \Bigg( \frac{5
\abs{\hat f}_{1, \tau}^2}{4 \nu \lambda} + \frac{5 C_1^2 e^{4} M_1^3
M_2}{4 \nu \lambda_1^{3/2}} \Bigg)^{1/2} <  \Bigg( \frac{\nu
d_2}{10} \Bigg)^{3/2}.
\end{equation*}
In order to satisfy the last inequality, we have to choose $\lambda$
large enough, such that
\begin{equation} \label{eq_lambda2}
\alpha^{2} \nu \lambda^{1/2} \lambda_1^{1/2} d_2
> \Bigg( \frac{C_4 \abs{\hat f}_{1, \tau}^2}{\nu \lambda} + \frac{C_5 M_1^3 M_2}{\nu \lambda_1^{3/2}}
\Bigg)^{1/3},
\end{equation}
for some absolute constants $C_4, C_5 > 0$. For such choice of
$\lambda$ we have $v(t)\in L^\infty([t_0, \infty), G^2_\tau)$, and
this proves the first part of the Theorem.

We are left to show that $v^\omega(x, t)$ is an asymptotic
approximation of the solution $u(x, t)$ of the NSV equation
(\ref{eq_nsv}). Let $z = u - v^\omega$, and denote $\bar z =
P_\lambda (u - v^{(2)})$, $\hat z = Q_\lambda u - \hat v$. Clearly,
by the construction and Theorem~\ref{thm_vm}, that
\begin{equation} \label{eq_zbar}
\lim_{t\to \infty} \nm{P_\lambda u(t) - \bar v(t)} = \lim_{t\to
\infty} \nm{\bar z(t)} = 0.
\end{equation}
Therefore, to prove (\ref{eq_asymgv}) we need to show that
\begin{equation*} \label{eq_zhat}
\lim_{t\to \infty} \nm{Q_\lambda u(t) - \hat v(t)} = \lim_{t\to
\infty} \nm{\hat z(t)} = 0.
\end{equation*}
Observe that $\hat z$ satisfies the equation
\[
{\hat z}_t + \nu A {\hat z} + \alpha^2 A {\hat z}_t + Q_\lambda
\big( B(u, z) + B(z, u) - B(z, z) \big) = 0, \;\;\;\;\; {\hat
z}(t_0) = Q_\lambda u(t_0).
\]
Taking an inner product of the last equation with $\hat z$ we get
\begin{equation} \label{zhat2}
\half \ddt \big( \abs{\hat z}^2 + \alpha^2 \nm{\hat z}^2 \big) + \nu
\nm{\hat z}^2 \le \bbabs{(B(\hat z, u), \hat z)} + \bbabs{ Q_\lambda
\big( B(u, \bar z) + B(\bar z, u) - B(z, \bar z), \hat z \big)}.
\end{equation}
The first summand on the right-hand side of equation (\ref{zhat2})
can be estimated as follows. Using (\ref{eq_nlin}),
Proposition~\ref{prop_poincare}, and relation (\ref{eq_v1bnd})
\[
\bbabs{(B(\hat z(t), u(t)), \hat z(t))} \le c \nm{u(t)} \abs{\hat
z(t)}^{1/2} \nm{\hat z(t)}^{3/2} \le \frac{c M_1}{\lambda^{1/4}
\lambda_1^{3/4}} \nm{\hat z(t)}^{2},
\]
for $t\ge t_0$. Plugging this inequality into equation
(\ref{zhat2}), and using relation (\ref{eq_nud2}), we get
\begin{equation} \label{zhat3}
\half \ddt \big( \abs{\hat z}^2 + \alpha^2 \nm{\hat z}^2 \big) + d_2
\Big( \nu - \frac{c M_1}{\lambda^{1/4} \lambda_1^{3/4}} \Big) \big(
\abs{\hat z}^2 + \alpha^2 \nm{\hat z}^2 \big) \le b(t),
\end{equation}
where
\[
b(t) = \abs{ ( B(u(t), \bar z(t)), \hat z(t))} + \abs{ (B(\bar z(t),
u(t)), \hat z(t))} + \abs{ (B(z(t), \bar z(t)), \hat z(t))}.
\]
Applying relation (\ref{eq_zbar}) and using the fact that $u$ is
bounded in the $V$ norm, we conclude that $b(t)\to 0$, as $t\to
\infty$. Therefore, applying Gronwall's Lemma~\ref{lem_1} to
equation (\ref{zhat3}) yields
\[
\lim_{t\to \infty} \nm{\hat z(t)} = 0,
\]
for $\lambda$ large enough, satisfying
\begin{equation} \label{eq_lamdetmodes}
\lambda > \lambda_1^{3} \Big( \frac{c M_1}{\nu} \Big)^4.
\end{equation}

Summarizing, the statement of the Theorem holds for $\lambda$ large
enough satisfying relations (\ref{eq_lambda1}), (\ref{eq_lambda2})
and (\ref{eq_lamdetmodes}).

\end{proof}

\section{Estimating the exponential decaying small scale}

As we have mentioned in the introduction, an additional goal of this
research is to provide further support for the proposal made in
\cite{CLT} that the NSV system (\ref{eq_nsv}), with the small
regularization parameter $\alpha$, can be used as a numerical model
for studying the original Navier-Stokes equations, and in particular
their statistical properties. Theorem~\ref{thm_gv} actually states
that the global attractor of the NSV system consists of real
analytic functions $u(x, t)$, whose Fourier spectrum $\hat u(k, t)$
satisfies the decay estimate
\[
\abs{\hat u(k, t)} \le c \abs{k}^{-2} e^{- \abs{k}/\lambda^{1/2} }.
\]
Therefore, following the ideas of \cite{DT95_Decay} (see also
\cite{HKR89_Smallscale}, \cite{HKR90_Smallscale} for a different
approach), the quantity $1 / \lambda^{1/2}$, can be naturally
identified as the \emph{exponential decaying length scale}, since
the exponential decay of the spectrum of $u$ is effective only at
high wavenumbers satisfying $\abs{k}
> \lambda^{1/2}$.

In the case of the Navier-Stokes equations the exponential decaying
length scale, and similarly the radius of analyticity of solutions,
can be identified with the \emph{smallest effective length scale} in
the turbulent flow (see, e.g., \cite{DT95_Decay}, \cite{FMRT01},
\cite{HKR89_Smallscale}, \cite{HKR90_Smallscale}). Classical
Kolmogorov theory of turbulence states that the smallest effective
length scale in the flow is proportional to
\[
\ell_K = \Big( \frac{\nu^3}{\epsilon} \Big)^{1/4},
\]
where
\[
\epsilon = \nu \ang{\nm{u}^2},
\]
is the mean energy dissipation rate, and $\ang{\cdot}$ denotes
either the long time average, or the ensemble average with respect
to the proper invariant probability measure. In \cite{DT95_Decay} it
was shown that for the solution $u(t)$ of the $3$D Navier-Stokes
equations, as long $\nm{u(t)}$ remains bounded uniformly on some
interval of time $[0, T]$, the smallest length scale of the
turbulent flow satisfies
\begin{equation} \label{eq_ellns}
\ell \sim L \Big( \frac{\ell_K}{L} \Big)^4,
\end{equation}
where in the definition of $\ell_K$, instead of the usual definition
of the energy dissipation rate $\epsilon$, the authors considered
the largest instantaneous energy dissipation rate on the time
interval $[t_1, T]$, on which the solution of the equations remains
regular
\[
\epsilon_{sup} = \sup_{t_1 \le t \le T} \nu \nm{v(t)}^2.
\]

In the case of the NSV system, similarly to the Navier-Stokes
equations, we can define $\ell_{NSV}$ -- the exponential decaying
length scale. In other words, $\ell_{NSV}$ is the largest length
scale below which an exponential decay of the spectrum of the
solutions of the NSV system lying on the global attractor becomes
effective. In this section we would like to derive a lower bound for
the $\ell_{NSV}$, similar to relation (\ref{eq_ellns}) for the $3$D
Navier-Stokes equations. For other estimates on a related smallest
length scale (via computation of the radius of analyticity of the
solutions) of the Navier-Stokes equations in $2$ and $3$ dimensions
see \cite{HKR89_Smallscale}, \cite{HKR90_Smallscale}, and
\cite{Ku99_Dissscale} (see also \cite{ITi07_Smallscales}). See also
\cite{CFMT85_Detmodes} and \cite{FMRT01} for other approach to this
subject.

The energy of the NSV system is defined as
\[
E(t) = \abs{u(t)}^2 + \alpha^2 \nm{u(t)}^2,
\]
which satisfies the balance
\[
\half \ddt E(t) = - \nu \nm{u}^2 + (f, u).
\]
We denote the mean rate of dissipation of energy for the NSV system
as
\[
\epsilon = \nu \ang{\nm{u}^2},
\]
where $\ang{\cdot}$ stands for the long-time average. Moreover, we
have the bound
\[
\epsilon \le \epsilon_{sup} := \nu M_1^2.
\]
In order to find a lower bound for the exponentially decaying length
scale of the NSV flow, we need to estimate the value of $\lambda$,
from the inequalities (\ref{eq_lambda1}), (\ref{eq_lambda2}) and
(\ref{eq_lamdetmodes}) in the proof of Theorem~\ref{thm_gv}, since
$\lambda^{-1/2}$ is a lower bound for the radius of analyticity of
the solutions of the NSV system lying on the attractor, and
therefore,
\[
\ell_{NSV} \ge \lambda^{-1/2}.
\]

First, note, that the condition (\ref{eq_lamdetmodes}) is satisfied
for
\begin{equation} \label{eq_lam1}
\lambda^{-1/2} \sim \frac{\nu^3}{L^3 \epsilon_{sup}} = L \Big(
\frac{\ell_K}{L} \Big)^4.
\end{equation}
Moreover, for a small viscosity $\nu$ and $\alpha$, we can estimate
$M_2$, using the expressions of Lemma~\ref{lem_2}, in the following
way
\[
M_2 \sim C_6 \frac{M_1^3}{\alpha^2 \nu^2 \lambda_1^{5/2}},
\]
where $C_6 > 0$ is an absolute constant. Therefore, the condition
(\ref{eq_lambda1}) holds if
\begin{equation} \label{eq_lam2}
\lambda^{-1/2} \sim \frac{\nu \lambda_1^{3/4}}{M_1^{1/2} M_2^{1/2}}
\sim \frac{\alpha \nu^2 \lambda_1^{2}}{M_1^{2}} \sim \frac{\alpha
\nu^3}{L^4 \epsilon_{sup}} = \alpha \Big( \frac{\ell_K}{L} \Big)^4,
\end{equation}
or, on the other hand if
\begin{equation} \label{eq_lam3}
\lambda^{-1/2} \sim \frac{\nu^2 \lambda_1^{3/2}}{M_1^2} \sim
\frac{\nu^3}{L^3 \epsilon_{sup}} = L \Big( \frac{\ell_K}{L} \Big)^4.
\end{equation}
Finally, we are left to check when the condition (\ref{eq_lambda2})
is satisfied. In order to do this, let us assume, as it is
conventionally done, that $\hat f = 0$, namely, $\lambda$ is chosen
large enough such that the forcing $f$ is supported on the modes
less than $\lambda^{-1}$. In addition, we assume that $\lambda >
\alpha^{-2}$, so that $d_2 \ge \half \alpha^{-2}$. In that case, the
condition (\ref{eq_lambda2}) becomes
\[
\lambda^{1/2} \sim \frac{M_1 M_2^{1/3}}{\nu^{4/3} \lambda_1} \sim
\frac{M_1^2}{\nu^{2} \alpha^{2/3} \lambda_1^{11/6}} \sim L^{11/3}
\alpha^{-2/3} \frac{M_1^2}{\nu^{2}},
\]
and we obtain the estimate
\begin{equation} \label{eq_lam4}
\lambda^{-1/2} \sim L^{1/3} \alpha^{2/3} \Big( \frac{\ell_K}{L}
\Big)^4.
\end{equation}

Combining relations (\ref{eq_lam1}), (\ref{eq_lam2}),
(\ref{eq_lam3}), and (\ref{eq_lam4}) we conclude that the
exponential decaying length scale of the NSV equations satisfies
\begin{equation} \label{eq_lnsv}
\ell_{NSV} \ge \min\{ L, \alpha, L^{1/3} \alpha^{2/3} \} \cdot \Big(
\frac{\ell_K}{L} \Big)^4.
\end{equation}
Note, that this estimate has the same asymptotic behavior as the
estimate of the characteristic length scale of the $3$D
Navier-Stokes equations obtained in \cite{DT95_Decay}, without
requiring any additional assumptions on the regularity of the flow
of the system (\ref{eq_nsv}).

\section{Radius of analyticity of stationary solutions}

At the end of the previous section we computed the exponential
decaying length scale of the NSV model by estimating the radius of
analyticity of the functions lying in the global attractor of the
system. A particular example of the functions lying on the attractor
are the stationary solutions of the system. The goal of this section
is to show that lower bounds for the exponential decaying length
scale of the stationary solutions of the NSV system are the same as
those obtained in the last section for the general element of the
global attractor. Observe that the NSV equations has the same
stationary solutions as the $3$D Navier-Stokes equations. All
calculation in this section are formal and can be rigorously
justified using the Galerkin approximation procedure. We are
following the ideas introduced in \cite{OliTi01_Stat}.

The steady state equation of (\ref{eq_nsv}) has the form
\begin{equation} \label{eq-stationary}
\nu A u + B(u, u) = f.
\end{equation}
Note the identity, where $\tau$ is now a dummy variable in the
interval $[0, \sigma]$
\begin{equation} \label{eq-0}
\frac{d}{d\tau} \nm{w}^2_{1, \tau} \equiv 2 \nm{w}^2_{3/2, \tau},
\end{equation}
for every time independent function $w\in G^{3/2}_\sigma$. Taking an
inner product of (\ref{eq-stationary}) with $A^{1/2} e^{2 \tau
A^{1/2}} u$, we obtain
\begin{align} \label{eq-1}
\nu \nm{u}^2_{3/2, \tau} & \le \bbabs{ (B(u, u), A^{1/2} u)_{0,
\tau} } + \bbabs{ (f, A^{1/2} u)_{0, \tau} }.
\end{align}

Let us assume that the forcing $f$ is supported on the first $N_f$
modes. Therefore, we can write
\[
\bbabs{ (f, A^{1/2} u)_{0, \tau} } \le e^{2\tau N_f^{1/2}} N_f^{1/2}
\bbabs{ (f, u) } \le e^{2\sigma N_f^{1/2}} N_f^{1/2} \nu \nm{u}^2,
\]
where the last inequality is the result of the energy conservation,
and we chose a large $\sigma$ to be determined later. Moreover, we
can estimate
\[
\bbabs{ (B(u, u), A^{1/2} u) }_{0, \tau} \le c \lambda_1^{3/4}
\nm{u}^2_{1, \tau} \nm{u}_{3/2, \tau} \le \frac{c^2}{2 \nu} \nm{ u
}^4_{1, \tau} + \frac{\nu}{2} \nm{ u }^2_{3/2, \tau}.
\]
Let us substitute the last two inequalities into equation
(\ref{eq-1}) to get
\begin{equation} \label{eq-2}
\half \nu \nm{u}^2_{3/2, \tau} \le \frac{c^2 \lambda_1^{3/2}}{2 \nu}
\nm{ u }^4_{1, \tau} + e^{2\sigma N_f^{1/2}} N_f^{1/2} \nu \nm{u}^2,
\end{equation}
which we can rewrite
\begin{equation} \label{eq-3}
\nm{u}^2_{3/2, \tau} \le \frac{c^2 \lambda_1^{3/2}}{\nu^2} \nm{ u
}^4_{1, \tau} + 2 e^{2\sigma N_f^{1/2}} N_f^{1/2} \nm{u}^2.
\end{equation}
Applying to the last inequality the identity (\ref{eq-0}) we obtain
\[
\frac{d}{d\tau} \nm{u}^2_{1, \tau} \le \frac{c^2
\lambda_1^{3/2}}{\nu^2} \nm{ u }^4_{1, \tau} + 2 e^{2\sigma
N_f^{1/2}} N_f^{1/2} \nm{u}^2,
\]
Denote
\[
y(\tau) = \frac{c^2 \lambda_1^{3/2}}{\nu^2} \nm{u}^2_{1, \tau},
\;\;\; F := \frac{c \lambda_1^{3/4}}{\nu} \sqrt{2} e^{\sigma
N_f^{1/4}} N_f^{1/2} \nm{u}.
\]
Therefore, $y(\tau)$ satisfies
\[
\dot y \le y^2 + F^2 \le (y + F )^2.
\]
Once again, denote, $z = y + F$, which satisfies
\[
\dot z \le z^2,
\]
therefore
\[
z(\tau) = y(\tau) + F \le \Big( z^{-1}(0) - \tau \Big)^{-1}.
\]
``Blow-up time'' is
\[
\tau_B > z^{-1}(0) = \frac{1}{\frac{c^2 \lambda_1^{3/2}}{\nu^2}
\nm{u}^2 + \frac{c \lambda_1^{3/4}}{\nu} \sqrt{2} e^{\sigma
N_f^{1/4}} N_f^{1/2} \nm{u} } \ge C \frac{\nu^2}{\lambda_1^{3/2}
\nm{u}^2} \ge L \Big( \frac{\ell}{L} \Big)^4,
\]
where we used the fact that $u$ is a steady state of the system
(\ref{eq_nsv}), and hence lies in the global attractor. As a result,
we get $u\in G^2_\tau$ for all $\tau < \tau_B$. Therefore, we showed
that the exponential decaying length scale of the stationary
solutions of the NSV system, satisfies the same bound
(\ref{eq_lnsv}) as the general solution of the NSV equations lying
on the global attractor. Moreover, this bound also holds for the
smallest length scales of the stationary solutions of the
Navier-Stokes equations -- similar to the general bound obtained in
\cite{DT95_Decay}.

\section{Conclusions}

We prove that the elements of the global attractor of the $3$D NSV
equations (\ref{eq_nsv}) with periodic boundary conditions, driven
by an analytic forcing, are analytic. A consequence of this result
is that the solutions of the 3D NSV system lying on the global
attractor have exponentially decaying spectrum, despite the fact
that the addition of the $-\alpha^2 \triangle u_t$ term changes the
parabolic character of the original Navier-Stokes equation, which
now starts to behave similar to a damped hyperbolic system.

An important consequence of our result is that the solutions of the
3D NSV system (\ref{eq_nsv}) lying on the global attractor posses a
dissipation range -- an exponentially decaying spectrum. This fact
provides an additional evidence that (\ref{eq_nsv}) with the small
regularization parameter $\alpha$ enjoys similar statistical
properties of the 3D Navier-Stokes equations, a fact that was first
suggested in \cite{CLT}.

Finally, following the ideas of \cite{OliTi01_Stat}, we have
computed a lower bound of the radius of analyticity of the steady
state solution of the NSV and Navier-Stokes equations. The bound
coincides with the one obtained for the general solutions of the
system (\ref{eq_nsv}) lying on the global attractor.

\section{Acknowledgements}

The work of V. K. Kalantarov was supported in part by The Scientific
and Research Council of Turkey, grant no. 106T337. B. Levant
acknowledges the hospitality of the Hausdorff Center for Mathematics
in Bonn University, where this work has started. The work of E. S.
Titi was supported in part by the NSF grants no.~DMS-0504619 and
no.~DMS-0708832, the ISF grant no.~120/6, and the BSF grant
no.~2004271.

\end{document}